\documentclass[11pt,reqno]{amsart}
\usepackage{xcolor}
\usepackage{graphicx}
\usepackage[T1]{fontenc}
\usepackage{amsmath,amssymb,amsthm}
\usepackage{esint}
\usepackage[a4paper,margin=2.8cm]{geometry}
\usepackage{tikz}
\usetikzlibrary{patterns}
\usepackage{subcaption}
\usepackage{mleftright} 
\mleftright

\usepackage[colorlinks=true,
            linkcolor=blue,
            citecolor=blue,
            urlcolor=blue]{hyperref}

\usepackage[backend=biber,
style=alphabetic,
sorting=nyt,
maxbibnames=10,
maxcitenames=10,
giveninits=true,
doi=false,
url=false]{biblatex}

\DeclareFieldFormat[article]{journaltitle}{#1}
\DeclareFieldFormat[book]{title}{#1}

\renewbibmacro{in:}{%
  \ifentrytype{article}{}{\printtext{\bibstring{in}\intitlepunct}}}

 \AtEveryBibitem{\clearfield{pubstate}}

\numberwithin{equation}{section}

\newcommand{\R}{\mathbb{R}}
\newcommand{\Z}{\mathbb{Z}}
\newcommand{\T}{\mathbb{T}}

\newcommand{\bw}{\left\lVert \rho \right\rVert_{B^{\omega,a}}}

\newcommand{\e}{{\rm e}}
\newcommand{\eps}{\varepsilon}

\newcommand{\divergence}{\operatorname{div}}

\theoremstyle{plain}
\newtheorem{thm}{Theorem}[section]
\newtheorem{cor}[thm]{Corollary}
\newtheorem{lem}[thm]{Lemma}
\newtheorem{prop}[thm]{Proposition}

\theoremstyle{definition}
\newtheorem{defn}[thm]{Definition}

\theoremstyle{remark}
\newtheorem{rem}[thm]{Remark}

\begin{document}

\numberwithin{equation}{section}

\title[Unbounded DiPerna-Lions commutators for Osgood-Sobolev vector fields]
{Unbounded DiPerna-Lions commutators for Osgood-Sobolev vector fields}

\author[L. Franchi]{Leonardo Franchi}
\address{Department of Pure Mathematics and Mathematical Statistics, Centre for Mathematical Sciences, Wilberforce Road, Cambridge CB3 0WA, United Kingdom}
\email{lf511@cam.ac.uk}
\author[C. J. P. Johansson]{Carl Johan Peter Johansson}
\address{Institute of Mathematics, EPFL, Station 8, 1015 Lausanne, Switzerland}
\email{carl.johansson@epfl.ch}

\begin{abstract}
In this paper we show that there exists a regime of Osgood-Sobolev vector fields in which the DiPerna-Lions scheme cannot be applied despite the recently proved well-posedness and propagation of regularity in \cite{dephilippis_franchi_2026_osgood}.
For Osgood-Sobolev vector fields we also show that the propagation of regularity proved in \cite{dephilippis_franchi_2026_osgood} is sharp and derive lower bounds on mixing.

\end{abstract}

\maketitle

\section{Introduction}
We study bounded weak solutions $\rho$ to the Cauchy problem for the
transport equation on the $d$-dimensional flat torus $\mathbb{T}^d$,
\begin{equation}\label{eq:TransportEquation}\tag{TE}
\left\{
\begin{array}{ll}
\partial_t \rho + b \cdot \nabla \rho = 0
& \text{in } (0,T)\times \mathbb{T}^d, \\
\rho(0,x)=\rho_0(x)
& x\in\mathbb{T}^d,
\end{array}
\right.
\end{equation}
where $b:[0,T]\times\mathbb{T}^d\to\mathbb{R}^d$ is a bounded,
divergence-free time-dependent vector field, and
$\rho_0 \colon \mathbb{T}^d \to \mathbb{R}$ is the initial datum.
In their seminal paper \cite{DiPernaLions1989}, DiPerna and Lions proved well-posedness of \eqref{eq:TransportEquation} in the class of bounded solutions for Sobolev vector fields. 
Later, Ambrosio \cite{Ambrosio2004} established the same result for vector fields of bounded variation.
A couple $(b,\rho)$ with $\operatorname{div}b = 0$ is said to be a weak solution to \eqref{eq:TransportEquation} if
\begin{equation}\label{eq:WeakTransport}
\int_0^T\int_{\mathbb{T}^d}
\rho\left[\partial_t\xi+b\cdot\nabla\xi\right]\,dx\,dt
=
-\int_{\mathbb{T}^d}\rho_0(x)\xi(0,x)\,dx, \quad \forall \xi\in C_c^\infty([0,T)\times\mathbb{T}^d).
\end{equation}

In this paper, we consider Osgood-Sobolev vector fields for which well-posedness and propagation of regularity for solutions to \eqref{eq:TransportEquation} have recently been established by the first author and De Philippis \cite{dephilippis_franchi_2026_osgood}.
Firstly, we investigate whether well-posedness can a posteriori be deduced from the classical DiPerna-Lions scheme leveraging the regularity of $\rho$ against the lack of regularity of $b$.
We show (see Theorem~\ref{thm:CommutatorBlowup}) that there exists a regime in which this is false, showing that \cite{dephilippis_franchi_2026_osgood} addresses a fundamental obstruction of the DiPerna-Lions scheme in the Osgood-Sobolev setting.
Secondly, we prove that the propagation of regularity proved in \cite{dephilippis_franchi_2026_osgood} is sharp (see Theorem~\ref{thm:loss-of-regularity}) and provide lower bounds on mixing in the setting of Osgood-Sobolev vector fields (see Theorem~\ref{thm:Mixing}).

\subsection{The DiPerna--Lions scheme}
The method introduced by DiPerna and Lions \cite{DiPernaLions1989} relies on renormalised solutions and commutator estimates. 
A solution $\rho$ is called renormalised if $\beta(\rho)$
is a distributional solution with initial datum $\beta(\rho_0)$ for every $\beta\in C^1(\mathbb{R})$,
and in the class of renormalised solutions \eqref{eq:TransportEquation} admits one unique solution. The remainder of the arguments consists therefore in proving that any bounded solution is renormalised. To prove this, we rely on commutator estimates.
Precisely, let $\varphi\in C_c^\infty(B_1)$ be a smooth radial mollifier with
$\int_{\mathbb{R}^d}\varphi=1$, set $\varphi_\eps(x):= \eps^{-d} \varphi\left({x}/{\eps}\right)$ and $\rho_\eps:= \rho \ast \varphi_\eps$.
The regularised solution satisfies
\[
\partial_t \rho_\eps+b\cdot\nabla \rho_\eps
=
-R_\eps(b,\rho),
\]
where the commutator is given by
\begin{equation}\label{eq:DiPernaLionsCommutator}
R_\eps(b,\rho)(x)
=
\int_{B_1}
\rho(x+\eps h)
\frac{b(x)-b(x+\eps h)}{\eps}
\cdot\nabla\varphi(h)\,dh.
\end{equation}
The chain rule then gives
\[
\partial_t\beta(\rho_\eps)
+b\cdot\nabla\beta(\rho_\eps)
=
-\beta'(\rho_\eps)R_\eps(b,\rho).
\]
Therefore, in order to prove that $\rho$ is renormalised, it is enough to show that 
$R_\eps(b,\rho) \to 0 \text{ in } L^1_{\mathrm{loc}}$ as $\eps \to 0$.
DiPerna and Lions proved this convergence when
\[
b\in L^\infty_{t,x}\cap L^1_t W^{1,p}_x
\qquad\text{for some }p\geq 1.
\]
Indeed, Sobolev regularity gives, for every fixed $h\in B_1$,
\[
\frac{b(x)-b(x+\eps h)}{\eps}
\longrightarrow
-\nabla b(x)h
\text{ in }L^1_{\mathrm{loc}}
\]
as $\eps\to0$. Passing to the limit in the expression for the commutator, the
resulting term is zero because $\operatorname{div}b=0$. Hence $ R_\eps(b,\rho) \to 0 \text{ in }L^1_{\mathrm{loc}}$, and therefore the solution is renormalised.
In this discussion, we assumed $b$ is incompressible for simplicity but the case of vector fields with bounded divergence can also be considered as is done in \cite{DiPernaLions1989}.

\subsection{The Osgood regime}
Osgood regularity is a relaxation of Lipschitz regularity under which the classical Cauchy-Lipschitz theory can be generalised \cite{Osgo98}. Precisely, a modulus of continuity $\omega \in C(\R_+;\R_+)$ is Osgood regular if
\begin{equation}
    \int_0^1 \frac{ds}{\omega(s)} = \infty
\end{equation}
and then $f \in C^\omega$ if $|f(x) - f(y)| \leq C \omega( \mathsf{d} (x,y))$ for all $x,y$, where $\mathsf{d}$ denotes the geodesic distance on $\T^d$. Recently, the conclusions of the DiPerna-Lions scheme, that is renormalisation and uniqueness of bounded solutions, have been extended to the setting of Osgood-Sobolev vector fields by the first author and De Philippis \cite{dephilippis_franchi_2026_osgood}.
A function $f$ is called Osgood-Sobolev regular if it satisfies 
\begin{equation}
    | f(x) - f(y) | \leq (g(x) + g(y) ) \omega ( \mathsf{d}(x,y) )
\end{equation}
for some $g \in L^p$. Then, it belongs to the so-called Haj\l{}asz space $\dot{M}^{\omega, p}$, precisely defined in Subsection~\ref{subsec:OsgSobSpaces}.
In \cite{dephilippis_franchi_2026_osgood}, vector fields $b \in L^1_t \dot{M}^{\omega, p}_x \cap L^\infty_{t,x}$ with $\omega$ Osgood and $p > 1$ (with $\omega$ satisfying a mild non-degeneracy condition, see \eqref{eq:LowerBoundOnMixingTechnicalAssumption}) are considered and this constitutes an extension of the DiPerna-Lions theory because with $\omega(s) = s$, we have $\dot{M}^{\omega, p} =  \dot{W}^{1,p}$.
The proofs in \cite{dephilippis_franchi_2026_osgood} are based on a
Littlewood-Paley decomposition and differ from the classical DiPerna-Lions scheme. In particular, they do not rely on proving that the DiPerna-Lions
commutator tends to zero in $L^1_{loc}$. Instead, the argument establishes the conservation of a
suitable weighted energy associated with the Littlewood-Paley components of
the solution, which yields uniqueness and renormalisation.
With the same approach, the first author and De Philippis also proved propagation of $B^{\omega, a}$ norms, defined in Subsection~\ref{subsec:RegMixNorms}.

\subsection{(Non-)Applicability of the DiPerna-Lions scheme in the Osgood setting} 
Given the new technique introduced in \cite{dephilippis_franchi_2026_osgood}, we believe it is natural to ask whether it addresses a fundamental obstruction of the DiPerna-Lions scheme in the Osgood-Sobolev setting or if the DiPerna-Lions scheme can be invoked a posteriori. At the end of this subsection, we will show that in some regime of Osgood moduli, the DiPerna-Lions scheme can indeed be invoked a posteriori. The mechanism behind this is that the lack of regularity of the vector field can be paired with the regularity propagated by the scalar. 
However, there is also a regime where this pairing is not sufficient and the DiPerna-Lions scheme cannot be invoked. The next theorem demonstrates this.

\begin{thm}\label{thm:CommutatorBlowup}
Let $1<p<\infty$, and let $\omega$ be an Osgood modulus of continuity. Let $G \colon [0,\infty) \to \R$ be defined by 
\begin{equation}
    G(r):=\int_{2^{-r}}^1\frac{ds}{\omega(s)}
\end{equation}
and assume that
\begin{equation}\label{eq:CommutatorGrowthHypothesis}
  \lim_{r\to\infty}\frac{2^r\omega(2^{-r})}{G(r)^{p-1}}=\infty.
\end{equation}
Then there exist
\[
  b\in L^\infty\big([0,2];
  \dot{M}^{\omega,p}(\T^2;\R^2)\cap L^\infty(\T^2;\R^2)\big), \quad \operatorname{div}b = 0,
\]
\[
  \rho\in C([0,2];L^1(\T^2))\cap L^\infty([0,2]\times\T^2),
  \qquad \rho(t,x)\in\{0,1\},
\]
solving \eqref{eq:TransportEquation} on $(0,2)\times\T^2$, with initial datum $\rho(0,\cdot)\in \operatorname{BV}(\T^2)$, such that
\begin{equation}\label{eq:CommutatorBlowupConclusion}
  \limsup_{\eps\downarrow0}\int_0^2
  \|R_\eps(b_t,\rho_t)\|_{L^1(\T^2)}\,dt=\infty.
\end{equation}
\end{thm}
 
Since the DiPerna-Lions scheme requires $\int_0^2 \| R_{\eps} (b_t, \rho_t) \|_{L^1} \, dt \to 0$ as $\eps \to 0$ and $(b,\rho)$ solve \eqref{eq:TransportEquation}, this shows non-applicability of the classical DiPerna-Lions scheme in this setting. 
Of course, \eqref{eq:CommutatorBlowupConclusion} does not exclude the existence of a subsequence of $\eps$ along which the sequence converges to $0$. 
Nevertheless, \eqref{eq:CommutatorBlowupConclusion} shows a clear obstruction to the classical DiPerna-Lions scheme and, to the best of our knowledge, is the first example of a pair consisting of a vector field and a renormalised solution to \eqref{eq:TransportEquation} whose commutators are not uniformly bounded in $\eps$.\\
We prove this theorem in Section~\ref{sec:NonConvergenceCommutators}.  
The vector field $b$ is constructed with infinitely many pieces, each consisting of two key parts. In $[0,1)$, each piece is a time-accelerated version of the mixing vector field by Alberti, Crippa and Mazzucato \cite{AlbertiCrippaMazzucato2019}, see Lemma~\ref{lemma:ConstructionOfmixingVectorField}.
In $[1,2]$, each piece is a space-periodised version of the vector field constructed in Lemma~\ref{lemma:SingleDiskCommutator}, making the commutator large.
\\
To clarify the meaning of Theorem~\ref{thm:CommutatorBlowup}, consider the simpler class of  moduli $\omega_\alpha$ defined by $\omega_\alpha(s) = ( \log({5}/{s}) )^\alpha s$ on $[0,1]$ with $\alpha \in [0,1]$.
Then, a calculation yields (for $r \gg 1$)
\begin{equation}\label{eq:AsymptoticsOfGForSimpleClassOfModuli}
    G(r) \sim r^{1 - \alpha} \text{ if $\alpha \in [0,1)$} \quad \text{and} \quad G(r) \sim \log (r) \text{ if $\alpha = 1$.}
\end{equation}
Thus, \eqref{eq:CommutatorGrowthHypothesis} holds whenever $\alpha > \frac{p-1}{p}$, showing the non-applicability of the DiPerna-Lions scheme in this regime.
We now provide a sketch showing that there exists a complementary regime in which the commutator does converge to zero. Assume $p \geq 2$ and recall the expression for $R_\eps(b,\rho)$ in \eqref{eq:DiPernaLionsCommutator} and observe that, since $\operatorname{div}b=0$, we have
\begin{equation}
\int_{B_1} \frac{b(t,x)-b(t,x+\eps h)}{\eps} \cdot\nabla\varphi(h)\,dh=0.
\end{equation}
Therefore
\begin{equation}
    R_\eps(b,\rho)(t,x) = \int_{B_1} \big( \rho(t,x+\eps h)-\rho(t,x) \big) \frac{b(t,x)-b(t,x+\eps h)}{\eps} \cdot \nabla\varphi(h)\,dh.
\end{equation}
In Appendix~\ref{sec:RegOfScalar}, we show that the regularity of $\rho$ is sufficient to deduce that
\begin{equation}
    |\rho(t,x)-\rho(t,y)| \lesssim_d \dfrac{ (q(t,x)+q(t,y))}{G(\log_2(1/ \mathsf{d}(x,y)))^{p-1}},
\end{equation}
with $q(t,\cdot)\in L^{\infty}_tL_x^{p^\prime}(\T^d)$.
Throughout this manuscript, we write $A\lesssim_\ast B$ if there exists a constant $C(\ast)>0$ depending on the parameters $\ast$ such that $A\le C(\ast) B$.
We write $A\sim_\ast B$ if $A\lesssim_\ast B$ and $B\lesssim_\ast A$.
At the same time, by assumption
\begin{equation}
    |b(t,x)-b(t,y)| \leq \omega( \mathsf{d}(x,y) )( g(t,x)+g(t,y) ).
\end{equation}
Combining the two last estimates and using H\"older's inequality, we obtain
\begin{equation}
    \|R_\eps(b,\rho)(t,\cdot)\|_{L^1} \lesssim_{d,p} \| b(t,\cdot) \|_{\dot{M}^{\omega, p}} \| q(t,\cdot) \|_{L^{p^\prime}} \frac{\omega(\eps)}{\eps G(\log_2(1/\eps))^{p-1}}.
\end{equation}
Thus, the DiPerna-Lions scheme can be applied if 
\begin{equation}
    \frac{\omega(\eps)}{\eps G(\log_2(1/\eps))^{p-1}} \to 0 \text{ as $\eps \to 0$.}
\end{equation}
Note that the quantity tending to zero is precisely the same as the one appearing in \eqref{eq:CommutatorGrowthHypothesis} after writing $\eps=2^{-r}$.
\subsection{Sharp regularity and lower bounds on mixing}
To complete the picture for vector fields satisfying the $L^p$-Osgood condition, we prove that the regularity result in \cite[Theorem 1.2]{dephilippis_franchi_2026_osgood} is sharp and establish a corresponding lower bound on mixing efficiency.
\begin{thm}\label{thm:loss-of-regularity}
Let $d \geq 2$, $1 < p < \infty$, and let $\omega$ be an Osgood modulus of continuity. 
Then there exist a bounded divergence-free vector field
\[
b \in L^\infty([0,\infty)\times\mathbb{T}^d;\mathbb{R}^d) \cap L^\infty\big([0,\infty);\dot{M}^{\omega,p}(\mathbb{T}^d;\mathbb{R}^d)\big)
\]
and a bounded weak solution $\rho \in L^\infty([0,\infty)\times\mathbb{T}^d)$ of \eqref{eq:TransportEquation} such that $\rho(0,\cdot) \in W^{1,d}(\T^d)$ but, for every $t > 0$ and every $a > \frac{p}{2}$, $\rho(t,\cdot) \notin B^{\omega,a}(\mathbb{T}^d)$.
\end{thm}

This proves that \cite[Theorem 1.2]{dephilippis_franchi_2026_osgood} is sharp.
The proof is based on the construction introduced in \cite{AlbertiCrippaMazzucato2019Loss}, but with each individual piece accelerated in time at small scales while maintaining $\dot{M}^{\omega, p}$ norms bounded.

Next, we show a lower bound on mixing efficiency, measured using negative Sobolev norms.
\begin{thm}\label{thm:Mixing}
Let $d \geq 2$, $1<p<\infty$, and let $\omega$ be an Osgood modulus of continuity such that there exist constants $C>0$ and $\beta>0$ satisfying
\begin{equation}\label{eq:LowerBoundOnMixingTechnicalAssumption}
    \omega(\lambda s) \leq C \lambda^\beta \omega(s) \quad \forall s,\lambda \in (0,1].
\end{equation}
Assume that $(b,\rho)$ solve \eqref{eq:TransportEquation} with $b$ divergence-free and $\rho$ mean-zero, with initial datum $0\neq\rho_0\in L^\infty(\mathbb{T}^d)\cap B^{\omega,(1+p)/4}(\mathbb{T}^d)$, and that
\begin{equation}
    b \in L^\infty([0,\infty) \times \mathbb{T}^d; \mathbb{R}^d) \cap L^\infty\bigl([0,\infty); \dot{M}^{\omega,p}(\mathbb{T}^d;\mathbb{R}^d)\bigr), \quad
    \rho \in L^\infty([0,\infty)\times\mathbb{T}^d).
\end{equation}
Then, for every $t\geq 0$,
\begin{equation}\label{eq:MixingLowerBoundInMixingThm}
\frac{\|\rho(t,\cdot)\|_{\dot H^{-1}}}{\|\rho_0\|_{L^2}} \geq \Phi_\omega\left( C\left( \frac{\|\rho_0\|_{B^{\omega,(1+p)/4}}}{\|\rho_0\|_{L^2}} \right)^{\frac{4}{1+p}} \right) \Phi_\omega\left( C\left( \frac{\|\rho_0\|_{L^\infty}}{\|\rho_0\|_{L^2}}\right)^{\frac{4}{1+p}} \int_0^t \|b(s,\cdot)\|_{\dot M^{\omega,p}}\,ds \right),
\end{equation}
where $\Phi_\omega(u):=e^{-G^{-1}(u)}$.
\end{thm}

The mixing rate $\Phi_\omega$ is sharp as shown in Remark~\ref{rem:SharpMixing}.
Equation~\eqref{eq:LowerBoundOnMixingTechnicalAssumption} is a technical assumption made in \cite{dephilippis_franchi_2026_osgood} satisfied by a large class of reasonable moduli of continuity. 
As for Theorem~\ref{thm:CommutatorBlowup}, to illustrate this result, we consider the simpler class of moduli $\omega_\alpha$ defined by $\omega_\alpha(s) = ( \log({5}/{s}) )^\alpha s$ on $[0,1]$ with $\alpha \in [0,1]$. Note that  \eqref{eq:LowerBoundOnMixingTechnicalAssumption} holds and recall \eqref{eq:AsymptoticsOfGForSimpleClassOfModuli} which after a calculation yields (for $t \gg 1$)
\begin{equation}
    G^{-1}(t) \sim t^{\frac{1}{1 - \alpha}} \text{ if $\alpha \in [0,1)$} \quad \text{and} \quad G^{-1}(t) \sim \e^t \text{ if $\alpha = 1$.}
\end{equation}
Hence, in the case $\alpha = 0$ where $\dot{M}^{\omega, p} = \dot{W}^{1,p}$, we recover the known fact that mixing occurs at most exponentially. On the other hand, in the log-Lipschitz setting, i.e. $\alpha = 1$, mixing occurs at most double-exponentially.

\subsection{Connections to previous literature}
Without aiming to be exhaustive, we outline some connections of this manuscript with the existing literature.
\begin{itemize}
    \item \textbf{Osgood regular vector fields.} In the same way that Sobolev vector fields are a relaxation of Lipschitz vector fields, Osgood-Sobolev vector fields are a relaxation of Osgood vector fields.
    The classical Osgood regularity \cite{Osgo98} guarantees existence and uniqueness of the associated ODE and holds in several important applications. Indeed, in the vorticity formulations of the two-dimensional Euler equations with initial bounded vorticity, the vector field is log-Lipschitz regular \cite{Yudovich} but fails to be uniformly Lipschitz in time as shown by the fact that the gradient of the vorticity can grow double-exponentially \cite{KiselevSverak}.
    Osgood vector fields also appear in various kinetic models \cite{CrippaInversiSaffirioStefani2024,JunneRege2025}. Further investigations into Osgood vector fields and their properties include \cite{AmbrBern08,LiLuo2015,ColoKuma26,JohaMesc26,BorrNari26}.
    
    \item \textbf{Regularity for rough vector fields.} In Theorem~\ref{thm:loss-of-regularity}, we proved that the regularity proved in \cite[Theorem 1.2]{dephilippis_franchi_2026_osgood} is sharp.
    In the classical Sobolev setting, the sharp regularity theory has been obtained in a series of works \cite{CrippaDeLellis2008, BrueNguyen2021, meyer2024propagation, HuysSaid26, BrueColoPhilJoha26} using a variety of methods. However, the case $p = 1$, that is the one concerned with $W^{1,1}$ vector fields, remains a challenging open problem. 
    \item \textbf{Lower bounds on mixing and Bressan's mixing conjecture.} In 2003, Bressan \cite{Bressan} conjectured a lower bound on mixing efficiency stating that mixing cannot occur faster than exponentially with exponential rate given by $\| \nabla b_t \|_{L^1}$. For a precise statement of the conjecture, we refer to \cite[Conjecture 1]{Bressan}. An $L^p$-version of the conjecture was proved in \cite{CrippaDeLellis2008}, see also \cite{IyerKiselevXu2014, Lege18}. The conjecture has also been verified for particular two-dimensional vector fields \cite{Cooperman, HadzSeegSmarStre18}. For the time-periodic vector fields, the conjecture has also been verified asymptotically \cite{BrueColoJoha25}.
\end{itemize}
\subsection*{Acknowledgments} 
L.F. acknowledges support from the Isaac Newton Trust through a Trinity Cambridge Research Studentship.
L.F. would like to thank Guido De Philippis for insightful discussions.
C.J. was supported by the Swiss State Secretariat for Education, Research and Innovation (SERI) under contract number MB22.00034 through the project TENSE.
C.J. would like to thank Maria Colombo and Elia Brué for their valuable guidance.
Part of this work was carried out while C.J. was visiting the Institute for Advanced Study, whose excellent working conditions and hospitality are gratefully acknowledged.
\subsection*{AI Disclosure Statement}
AI tools were used for copyediting, proofreading, simplifying the original proof of Lemma~\ref{lemma:SingleDiskCommutator}, and generating Figure~\ref{fig:pinching}. 

\section{Preliminaries}
\subsection{Regularity and mixing norms}\label{subsec:RegMixNorms}
The relevant norms to measure regularity and mixing are defined using Littlewood-Paley blocks. For a detailed treatment of this topic, we refer to \cite{grafakos2014classical, BahouriCheminDanchin2011}. Any function $f \in L^2$ can be written, up to a zero measure set, as
\begin{equation}\label{eq:LittlewoodPaleyDecomposition}
    f = f_0 + \sum_{k = 1}^{\infty} f_k \quad \text{with } \operatorname{supp}( \widehat{f_k} ) \subseteq B_{2^{k}} \setminus B_{2^{k-2}} \text{ if $k \geq 1$ and $\operatorname{supp}( \widehat{f_0} ) \subseteq B_1$},
\end{equation}
where $\widehat{f_k}$ denotes the Fourier transform of $f_k$.
The functions $f_k$ are called the Littlewood--Paley blocks of $f$.
For $k \geq 1$, we also introduce the notation
\begin{equation*}
    f_k^{\geq} := \sum_{j=k}^{\infty} f_j, \qquad f_k^{\leq} := f - f_{k+1}^{\geq}.
\end{equation*}
Moreover, the decomposition in \eqref{eq:LittlewoodPaleyDecomposition} is made so that
\begin{equation}\label{eq:PlancheralEquivalence}
    \| f \|_{L^2}^2 \sim_d \sum_{k = 0}^\infty \| f_k \|_{L^2}^2.
\end{equation}
Then, the regularity and mixing semi-norms can be defined as
\begin{equation*}
    \| f \|_{B^{\omega, a} }^2 := \sum_{k = 0}^\infty G(k)^{2 a} \| f_k \|_{L^2}^2 \qquad \| f \|_{\dot{H}^{-1} }^2 := \sum_{k = 1}^\infty 2^{- 2 k} \| f_k \|_{L^2}^2.
\end{equation*}

\subsection{Osgood-Sobolev spaces}\label{subsec:OsgSobSpaces}

\begin{defn}
Let $\omega:[0,\infty)\to[0,\infty)$ be a modulus of continuity, i.e., a continuous, increasing and concave function satisfying
$\omega(0)=0$ and $\lim_{r \to \infty} \omega(r) = \infty$. A measurable function $f \colon \T^d \to \R$ is said to belong to the Haj\l{}asz space $\dot{M}^{\omega, p}(\T^d)$ if there exists a function $g \in L^p(\T^d)$ such that
    \begin{equation}\label{eq:HajlaszDefinitionIneq}
        |f(x) - f(y)| \leq ( g(x) + g(y) ) \,\omega ( \mathsf{d}(x,y) ) \quad \forall x,y \in \T^d.
    \end{equation}
    For such functions, we define the semi-norm
    \begin{equation}
        \|f\|_{\dot{M}^{\omega,p}} := \inf \left\{ \|g\|_{L^p(\mathbb{T}^d)} :  \text{ $g$ satisfies \eqref{eq:HajlaszDefinitionIneq} } \right\}.
    \end{equation}
\end{defn}

The following lemma shows that $\dot{M}^{\omega,p}$
interpolates between $L^p$ and $W^{1,p}$.
\begin{lem}\label{lemma:InterpolationOsgood}
Let $1<p<\infty$. For every non-constant $f \in W^{1,p}(\T^d)$, we have
\begin{equation}
\label{eq:InterpolationEstimate}
\| f \|_{\dot{M}^{\omega,p}} \leq C(d,p) \frac{{\|f\|_{L^p}}}{\omega\left(\frac{\|f\|_{L^p}}{\|\nabla f\|_{L^p}}\right)}.
\end{equation}
\end{lem}

\begin{proof}
By the classical Haj\l{}asz inequality (see \cite{Hajlasz1996}),
for almost every $x,y\in\mathbb{T}^d$, we have
\begin{equation}\label{eq:interpolation-proof:hajlasz}
  |f(x)-f(y)| \lesssim_d \big (M|\nabla f|(x)+M|\nabla f|(y) \big)\mathsf{d}(x,y),
\end{equation}
where $M$ denotes the Hardy-Littlewood maximal operator on $\mathbb{T}^d$.
We also have the elementary estimate
\begin{equation}\label{eq:interpolation-proof:triangle}
  |f(x)-f(y)| \leq |f(x)|+|f(y)|.
\end{equation}

We now fix a parameter $r>0$. Concavity of $\omega$ and $\omega(0) = 0$ give
$\omega(s) \geq \frac{s}{r} \omega(r)$ for all $0\le s\le r$.
Thus, when $0<\mathsf{d}(x,y)<r$,
\eqref{eq:interpolation-proof:hajlasz} gives
\[
  |f(x)-f(y)| \lesssim_d \frac{r}{\omega(r)} \big(M|\nabla f|(x)+M|\nabla f| (y)\big) \omega ( \mathsf{d}(x,y) ).
\]
On the other hand, when $\mathsf{d}(x,y)\ge r$, monotonicity of $\omega$ and
\eqref{eq:interpolation-proof:triangle} imply
\[
  |f(x)-f(y)| \leq \frac{|f(x)|+|f(y)|}{\omega(r)} \omega ( \mathsf{d}(x,y) ).
\]
The non-negative function
\[
  g_r(x):=\frac{|f(x)|+C(d)r\,M|\nabla f|(x)}{\omega(r)}
\]
satisfies \eqref{eq:HajlaszDefinitionIneq} and, since $1<p<\infty$, the $L^p$-boundedness of $M$ gives
\begin{equation}
\label{eq:InterpolationEstimateAtScale}
  \|f\|_{\dot{M}^{\omega,p}} \leq \|g_r\|_{L^p} \lesssim_{d,p} \frac{\max\{\|f\|_{L^p},r\|\nabla f\|_{L^p}\}}{\omega(r)}.
\end{equation}
Finally, since $f$ is nonconstant, we can choose $r:=\frac{\|f\|_{L^p}}{\|\nabla f\|_{L^p}}$ in \eqref{eq:InterpolationEstimateAtScale} to get \eqref{eq:InterpolationEstimate}.
\end{proof}
We will also need a localization lemma for functions in $\dot{M}^{\omega, p}$.
\begin{lem}\label{lem:DisjointnessOfSupportsForHajSpaces}
    Let $1 < p < \infty$, $f \in \dot{M}^{\omega, p} \cap L^p (\T^d)$ with $\operatorname{supp}(f) \subseteq B_r(0)$ with $r > 0$. Then there exists $h \in L^p(\T^d)$ satisfying \eqref{eq:HajlaszDefinitionIneq} and such that 
    \begin{equation}
        \| h \|_{L^p} \leq \| f \|_{\dot{M}^{\omega,p}} + \frac{\| f \|_{L^{p}}}{\omega(r)} \quad \text{and} \quad \operatorname{supp}(h) \subseteq B_{2r}(0).
    \end{equation}
\end{lem}
\begin{proof}
    Let $g \in L^p(\T^d)$ be any function satisfying \eqref{eq:HajlaszDefinitionIneq} such that $\| g \|_{L^p} = \| f \|_{\dot{M}^{\omega,p}}$ and set $h = g 1_{B_{2r}(0)} + \frac{|f|}{\omega(r)}$.
    We will prove that
    \begin{equation}\label{eq:WhatToProveInCompactSupportLemma}
        |f(x) - f(y)| \leq ( h(x) + h(y) ) \omega( \mathsf{d}(x,y) ) \quad \forall x,y \in \T^d.
    \end{equation}
    If $x,y \in B_{2r}(0)$ or $x,y \not\in B_{2r}(0)$, then \eqref{eq:WhatToProveInCompactSupportLemma} trivially holds.
    It only remains to consider the mixed case where $x \in B_{2r}(0)$ and $y \not\in B_{2r}(0)$. 
    If in that case, $x \not\in B_r(0)$, then \eqref{eq:WhatToProveInCompactSupportLemma} holds and therefore we assume $x \in B_r(0)$. 
    Then
    \begin{equation}
        |f(x) - f(y)| = |f(x)| \leq \frac{|f(x)|}{\omega(r)} \omega( \mathsf{d}(x,y) ) \leq h(x) \omega( \mathsf{d}(x,y) ).
    \end{equation}
\end{proof}

\subsection{Mixing vector fields}
We work on $\R^2$ with compactly supported vector fields; we will restrict to the torus later. Throughout this subsection, unspecified spatial domains are $\R^2$, and \eqref{eq:TransportEquation} is understood on $\R^2$.

For every integer $m\geq0$, let $Q_m = \left(-2^{-(m+1)},2^{-(m+1)}\right)^2$.
Thus $Q_0$ is the unit square and $Q_m$ has side length $2^{-m}$.
The following elementary mixing step is the pinching construction of Alberti, Crippa and Mazzucato \cite[Section 5]{AlbertiCrippaMazzucato2019}.

\begin{lem}\label{lemma:ACM-Lemma}
There exist a bounded divergence-free vector field $U \colon [0,1] \times \R^2 \to \R^2$, vanishing outside $Q_0$, and open sets $F(t)\subseteq Q_0$, $0\leq t\leq 1$, with the following properties:
\begin{enumerate}
    \item For every $1<p<\infty$, $U|_{Q_0}\in L^\infty([0,1];W^{1,p}_0(Q_0;\R^2))$, and
    \[
      \| U \|_{L^\infty([0,1];L^\infty)} \leq C, \qquad \| U \|_{L^\infty([0,1];W^{1,p})} \leq C(p).
    \]
    \item The map $t \mapsto \mathbf{1}_{F(t)}$ belongs to $C([0,1];L^1)$ and solves
    \[
      \partial_t\mathbf{1}_{F(t)} + \divergence \big( U\mathbf{1}_{F(t)} \big) = 0
    \]
    in the sense of distributions on $(0,1)\times\R^2$.
    \item For every $t\in[0,1]$, $\mathcal{L}^2(F(t))= \mathcal{L}^2(F(0))$, and
    \[
     F(0)=B_{1/4}(0), \qquad F(1) = \bigcup_{\sigma\in\{-1,1\}^2} B_{1/8}\left(\frac{\sigma}{4}\right).
    \]
\end{enumerate}
\end{lem}

For integers $m,j \geq 0$, let $\mathcal{Q}_{m,j}$ be the collection of $4^j$ open dyadic squares of $Q_m$, each of side length $2^{-(m+j)}$.
For each $S\in\mathcal{Q}_{m,j}$, write $x_S$ for its centre, and set

\begin{equation}
E_{m,j}=\bigcup_{S\in\mathcal{Q}_{m,j}} B_{2^{-(m+j+2)}}(x_S).
\end{equation}

\newcommand{\pinchingpanel}[1]{%
  \begin{tikzpicture}[
    x=3.6cm, y=3.6cm,
    boundary/.style={draw=black, line width=0.6pt},
    disk/.style={draw=black, line width=0.45pt,
      pattern=north east lines, pattern color=black!65}
  ]
    \draw[boundary] (-0.5,-0.5) rectangle (0.5,0.5);
    \pgfmathtruncatemacro{\n}{2^(#1)}
    \pgfmathsetmacro{\radius}{1/(4*\n)}
    \foreach \i in {1,...,\n} {
      \foreach \k in {1,...,\n} {
        \draw[disk]
          ({-0.5+(\i-0.5)/\n},{-0.5+(\k-0.5)/\n})
          circle[radius=\radius];
      }
    }
  \end{tikzpicture}%
}

\begin{figure}[htbp]
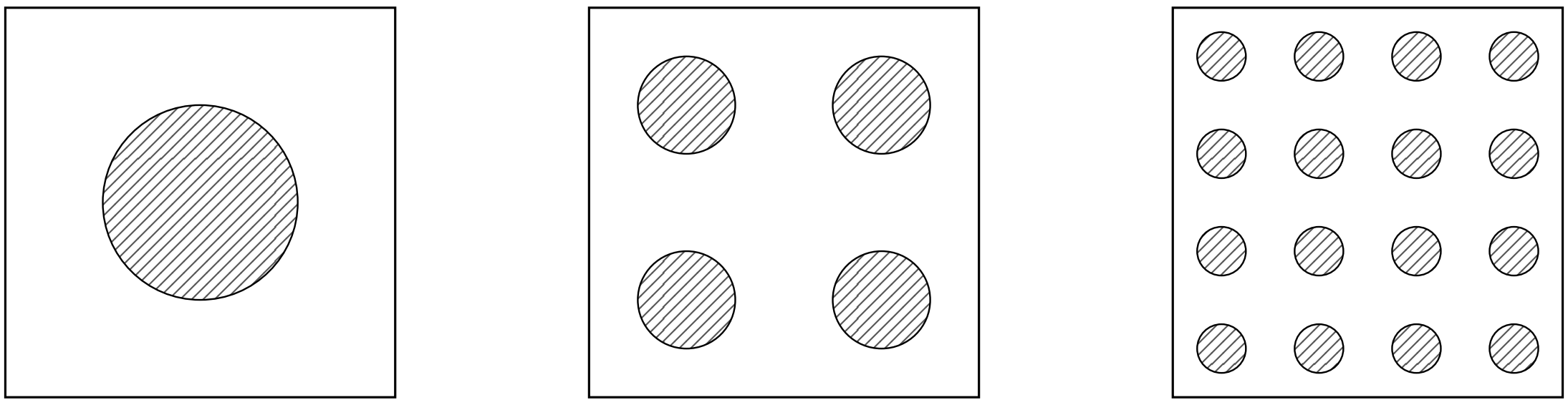

  \centering
  \begin{subfigure}[t]{0.3\textwidth}
    \centering
    \pinchingpanel{0}
    \caption{The set $E_{m,0}$}
    \label{fig:pinching-initial}
  \end{subfigure}\hfill
  \begin{subfigure}[t]{0.3\textwidth}
    \centering
    \pinchingpanel{1}
    \caption{The set $E_{m,1}$}
    \label{fig:pinching-first}
  \end{subfigure}\hfill
  \begin{subfigure}[t]{0.3\textwidth}
    \centering
    \pinchingpanel{2}
    \caption{The set $E_{m,2}$}
    \label{fig:pinching-second}
  \end{subfigure}
  \caption{The iterated construction in $Q_m$.
    At each step, the construction replaces the disk in every dyadic
    sub-square by four disks of half the radius. 
    The area of the shaded set is preserved. For related illustrations, see \cite[Figures 1 and 3]{AlbertiCrippaMazzucato2019}.}
  \label{fig:pinching}
\end{figure}

\begin{lem}\label{lemma:ConstructionOfmixingVectorField}
There exists a universal constant $c_0 \in (0,1)$ with the following property. Let $1 < p < \infty$, $\beta \in (0,1)$ and let $m,J \geq 1$ be integers. 
If $2^{-m} \leq c_0 \beta$ and 
 \begin{equation}\label{eq:ConditionOnJ}
   2^{-2m/p} \big( G(m+J)-G(m) \big) \leq c_0 \beta,
 \end{equation}
then there exist a divergence-free vector field $u \colon [0,1] \times \R^2 \to \R^2$, vanishing outside $Q_m$, and open sets $E(t)\subseteq Q_m$ such that:
\begin{enumerate}
  \item\label{item:ConstructionOfmixingVectorFieldItem1}
  $u|_{Q_m}\in L^\infty([0,1];W^{1,p}_0(Q_m;\R^2))$, and
  \begin{equation}
  \|u\|_{L^\infty([0,1];L^\infty)} \leq C\beta, \qquad \|u\|_{L^\infty([0,1];\dot{M}^{\omega,p})} \leq C(p)\beta.
  \end{equation}
  \item\label{item:ConstructionOfmixingVectorFieldItem2}
  The map $t\mapsto\mathbf{1}_{E(t)}$ belongs to $C([0,1];L^1)$ and solves \eqref{eq:TransportEquation} with $u$, i.e.
  \[
    \partial_t \mathbf{1}_{E(t)} + u \cdot \nabla \mathbf{1}_{E(t)} =0.
  \]
  Moreover, $\mathcal{L}^2(E(t))=\mathcal{L}^2(E(0))$ for every $t\in[0,1]$, and
  \[
    E(0)=B_{2^{-(m+2)}}(0),\qquad E(1)=E_{m,J}.
  \]
\end{enumerate}
\end{lem}

\begin{proof}
Set $a_j=2^{-(m+j)}$ for $0\leq j\leq J$. 
Let $\tau_0,\ldots,\tau_{J-1}>0$ be durations to be chosen, and define $t_j=\sum_{k=0}^{j-1}\tau_k$ for all $0 \leq j \leq J$.
Set $I_j=[t_j,t_{j+1})$ for $0\leq j<J$.
We will choose the durations so that $t_J\leq1$.
Let $U$ and $F$ be given by Lemma~\ref{lemma:ACM-Lemma}.
For $0\leq j<J$, $t\in I_j$ and $x\in\R^2$, define
\begin{equation}\label{eq:MixingVelocityDefinition}
  u(t,x)=\frac{a_j}{\tau_j} \sum_{S\in\mathcal{Q}_{m,j}} U \left( \frac{t-t_j}{\tau_j}, \frac{x-x_S}{a_j} \right).
\end{equation}
Set $u(t,\cdot) = 0$ for $t_J \leq t \leq1$.
Each term is a globally defined divergence-free field supported in its own square and therefore $u|_{Q_m}\in L^\infty([0,t_J];W^{1,p}_0(Q_m;\R^2))$.
Disjointness of the squares and the identity $4^j a_j^2=2^{-2m}$ give, for $t\in I_j$,
\begin{equation}\label{eq:MixingLpBound}
  \|u(t)\|_{L^p}\leq C\frac{a_j\,2^{-2m/p}}{\tau_j}, \quad
  \|\nabla u(t)\|_{L^p}\leq C(p)\frac{2^{-2m/p}}{\tau_j}, \quad
  \|u(t)\|_{L^\infty}\leq C\frac{a_j}{\tau_j}.
\end{equation}
Applying \eqref{eq:InterpolationEstimateAtScale} with $r = a_j$ and using \eqref{eq:MixingLpBound}, we obtain
\begin{equation}\label{eq:OsgoodBoundOnMixingVectorField}
  \|u(t)\|_{\dot{M}^{\omega,p}} \leq C(p)\frac{a_j\,2^{-2m/p}}{\tau_j\omega(a_j)}.
\end{equation}
To control both the $L^\infty$ and $\dot{M}^{\omega,p}$ norms, choose
\[
  \tau_j=\frac{a_j}{\beta} \left(1+\frac{2^{-2m/p}}{\omega(a_j)}\right), \qquad 0\leq j<J.
\]
Then \eqref{item:ConstructionOfmixingVectorFieldItem1} follows.
We now check that $t_J\leq1$. 
Monotonicity of $\omega$ gives $\frac{a_j}{\omega(a_j)} \leq 2 ( G(m + j + 1) - G(m + j) )$.
Since $\sum_{j=0}^{J-1}a_j\leq2^{1-m}$, it follows that
\begin{equation}\label{eq:MixingTotalTime}
  t_J\leq\frac{2}{\beta} \left( 2^{-m}+2^{-2m/p} \big( G(m+J)-G(m) \big) \right)\leq1,
\end{equation}
provided $c_0=1/4$ in \eqref{eq:ConditionOnJ} and $2^{-m}\leq\frac{\beta}{4}$.
We now define $E(t)$ as
\begin{equation}\label{eq:MixingSetsDefinition}
  E(t) = \bigcup_{S\in\mathcal{Q}_{m,j}} \left(x_S+a_j F\left(\frac{t-t_j}{\tau_j}\right)\right) \quad \forall t\in[t_j,t_{j+1}].
\end{equation}

At $t=t_j$ this is $E_{m,j}$. 
At $t=t_{j+1}$, the four disks in each square are precisely the disks associated with its four dyadic children, so the result is $E_{m,j+1}$. 
Thus the definitions agree at all $t_j$. Set $E(t)=E_{m,J}$ for $t_J\leq t\leq1$.
It follows from Lemma~\ref{lemma:ACM-Lemma} that $\mathbf{1}_{E(\cdot)}\in C([0,1];L^1)$ and $\mathcal{L}^2(E(t))=\mathcal{L}^2(E(0))$ for every $t\in[0,1]$.
Rescaling the equation in Lemma~\ref{lemma:ACM-Lemma} and using $L^1$ continuity at the joining times yields the transport equation on $(0,1)\times\R^2$, completing \eqref{item:ConstructionOfmixingVectorFieldItem2}.
\end{proof}

In the classical Sobolev setting, the vector fields constructed here are vector fields for which there is one solution to the transport equation that mixes exponentially and it corresponds to the construction in \cite{AlbertiCrippaMazzucato2019}. The authors also present a construction of a smooth vector field with smooth solution that mixes exponentially (however higher order norms of the vector field grow to infinity as time goes to infinity). Precisely they prove the following: 

\begin{thm}[{\cite[Theorem 6]{AlbertiCrippaMazzucato2019Loss}}]\label{lipschitz}
Let $d\ge 2$ and let $Q\subset\mathbb R^d$ be a fixed cube. There exist a smooth divergence-free velocity field $v$ and a smooth nontrivial solution $\rho$ of the transport equation, both spatially supported in $Q$, such that
\[
\rho_0:=\rho(0,\cdot)\in C_c^\infty(Q),
\]
with $\rho_0 \not \equiv 0$ and a mean-zero function.
Moreover, for every $t \geq 0$ and $x,y \in \R^d$, we have
\begin{align}
    \| v(t,\cdot) \|_{\operatorname{Lip}} &\leq C \label{eq:LipRegOfSmoothBuildingBlock} \\
    \|v(t, \cdot)\|_{L^\infty} +\|\rho(t, \cdot)\|_{\dot H^{-1}} \label{eq:DecayOfSmoothBuildingBlock} &\leq Ce^{-t}.
\end{align}
\end{thm}

We will use rescaled and accelerated versions of these vector fields in the proof of Theorem~\ref{thm:loss-of-regularity}.
\section{Non-convergence of commutators}\label{sec:NonConvergenceCommutators}

\subsection{Building blocks}

\begin{lem}\label{lemma:SingleDiskCommutator}
 There exist constants $\eps_0\in(0,1/5)$ and $c_0>0$, and a divergence-free vector field $W\in C_c^\infty(B_{1/3}(0)\setminus\overline{B_{1/4}(0)};\R^2)$ such that
 \begin{equation}\label{eq:CommutatorLowerBoundInLemma}
  \|R_{\eps_0}(W,\mathbf{1}_{B_{1/4}(0)})\|_{L^1} \geq c_0.
 \end{equation}
\end{lem}

\begin{proof}
For any $W$ vanishing on $B_{1/4}(0)$, the commutator definition gives
\[
  R_{\eps_0}(W,\mathbf{1}_{B_{1/4}(0)})(x) = -W(x) \cdot \nabla \big(\varphi_{\eps_0}*\mathbf{1}_{B_{1/4}(0)}\big)(x).
\]
If $\varphi$ is radial, the gradient is radial, so it suffices to choose a divergence-free $W$ with a nonzero radial component where this gradient does not vanish. For a general mollifier, we proceed as follows.

For $\eps_0>0$ sufficiently small, the convolution $\varphi_{\eps_0} \ast \mathbf{1}_{B_{1/4}(0)}$ is larger than $1/4$ somewhere on $\partial B_{1/4}(0)$ and vanishes outside $B_{1/3}(0)$. 
Thus, there exists $x_0 \in B_{1/3}(0) \setminus \overline{B_{1/4}(0)}$ with $\nabla \big(\varphi_{\eps_0}*\mathbf{1}_{B_{1/4}(0)}\big)(x_0) \in \R^2 \setminus B_1(0)$.
By continuity, there exists $\delta_0=\delta_0(\eps_0) > 0$ such that $\overline{B_{\delta_0}(x_0)} \subset B_{1/3}(0) \setminus \overline{B_{1/4}(0)}$ and
\[
  \left|\nabla\bigl(\varphi_{\eps_0}*\mathbf{1}_{B_{1/4}(0)}\bigr)(x)- \nabla\bigl(\varphi_{\eps_0}*\mathbf{1}_{B_{1/4}(0)}\bigr)(x_0) \right| \leq \eps_0 \qquad \text{for }x\in B_{\delta_0}(x_0).
\]
Selecting $W \in C_c^\infty(B_{1/3}(0)\setminus\overline{B_{1/4}(0)};\R^2)$ incompressible and such that $W = \nabla\bigl(\varphi_{\eps_0}*\mathbf{1}_{B_{1/4}(0)}\bigr)(x_0)$ on $B_{\delta_0 / 2}(x_0)$ implies \eqref{eq:CommutatorLowerBoundInLemma}.
\end{proof}

\begin{cor}\label{cor:MultipleRotationsInASquare}
For any integers $m\geq0$ and $J\geq1$, there exists a divergence-free vector field $W_{m,J}\in C_c^\infty(Q_m;\R^2)$ such that:
\begin{enumerate}
  \item\label{item:MultipleRotationsInASquareItem1}
  For every $1<p<\infty$,
  \begin{equation}\label{eq:MultipleRotationsNormBounds}
    \|W_{m,J}\|_{L^p}\leq C2^{-(m+J)-2m/p}, \quad
    \|\nabla W_{m,J}\|_{L^p}\leq C2^{-2m/p}, \quad
    \|W_{m,J}\|_{L^\infty}\leq C2^{-(m+J)},
  \end{equation}
  and
  \begin{equation}\label{eq:MultipleRotationsOsgoodBound}
    \|W_{m,J}\|_{\dot{M}^{\omega,p}} \leq C(p)\frac{2^{-(m+J)-2m/p}}{\omega(2^{-(m+J)})};
  \end{equation}
  \item\label{item:MultipleRotationsInASquareItem2}
  $\mathbf{1}_{E_{m,J}}$ is a stationary solution of \eqref{eq:TransportEquation} with $W_{m,J}$;
  \item\label{item:MultipleRotationsInASquareItem3}
  With $\eps_0,c_0$ as in Lemma~\ref{lemma:SingleDiskCommutator},
  \begin{equation}\label{eq:MultipleRotationsCommutatorEstimate}
    \|R_{2^{-(m+J)}\eps_0}(W_{m,J},\mathbf{1}_{E_{m,J}})\|_{L^1} \geq c_0 2^{-2m}.
  \end{equation}
\end{enumerate}
\end{cor}

\begin{proof}
Let $W$ be given by Lemma~\ref{lemma:SingleDiskCommutator}.
For $x\in\R^2$, define
\begin{equation}\label{eq:DefinitionOfWEllJ}
  W_{m,J}(x)=2^{-(m+J)}\sum_{S\in\mathcal{Q}_{m,J}} W\left( \frac{x-x_S}{2^{-(m+J)}} \right).
\end{equation}
The field is smooth, divergence-free and compactly supported in $Q_m$.
Disjointness of the cell supports and $4^J 2^{-2(m+J)} = 2^{-2m}$ give \eqref{eq:MultipleRotationsNormBounds}.
These estimates and Lemma~\ref{lemma:InterpolationOsgood} yield \eqref{eq:MultipleRotationsOsgoodBound}.
Since $W_{m,J}$ vanishes on $E_{m,J}$, $\mathbf{1}_{E_{m,J}}$ is a stationary solution of the transport equation.
Using $W_{m,J}\mathbf{1}_{E_{m,J}}=0$ and rescaling the commutator gives
\begin{equation}\label{eq:CommutatorOfWEllJ}
  R_{2^{-(m+J)}\eps_0}(W_{m,J},\mathbf{1}_{E_{m,J}})(x) = \sum_{S\in\mathcal{Q}_{m,J}} R_{\eps_0}(W,\mathbf{1}_{B_{1/4}(0)}) \left(\frac{x-x_S}{2^{-(m+J)}}\right).
\end{equation}
Since the terms in the sum have disjoint supports, \eqref{eq:MultipleRotationsCommutatorEstimate} follows from Lemma~\ref{lemma:SingleDiskCommutator} and $\#\mathcal{Q}_{m,J} = 4^J$.
\end{proof}

\begin{prop}\label{prop:CombinedBuildingBlock}
Let $1<p<\infty$, $\beta\in(0,1)$, and let $m,J\geq1$ be integers satisfying the hypotheses of Lemma~\ref{lemma:ConstructionOfmixingVectorField}.
There exist a divergence-free vector field $b_m \colon [0,2] \times \R^2 \to \R^2$ and $\rho_m \colon [0,2] \times \R^2 \to \{0,1\}$, both vanishing outside $Q_m$, such that:
\begin{enumerate}
  \item
  \begin{equation}\label{eq:CombinedBlockNormBounds}
    \|b_m\|_{L^\infty([0,2];L^\infty)}\leq C\beta,
    \qquad
    \|b_m\|_{L^\infty([0,2];\dot{M}^{\omega,p})}\leq C(p)\beta;
  \end{equation}
  \item
  $\rho_m\in C([0,2];L^1)$ solves \eqref{eq:TransportEquation} with $b_m$ and initial datum $\rho_m(0,\cdot)=\mathbf{1}_{B_{2^{-(m+2)}}(0)}$.
  Moreover, with $\eps_0,c_0$ as in Lemma~\ref{lemma:SingleDiskCommutator}, for every $t\in(1,2]$,
  \begin{equation}\label{eq:OneFullBlockCommutatorEstimate}
    \left\|R_{2^{-(m+J)}\eps_0}\big(b_m(t,\cdot),\rho_m(t,\cdot) \big)\right\|_{L^1}
    \geq c_0\beta\frac{2^{J-m}\omega(2^{-(m+J)})}{\omega(2^{-(m+J)})+2^{-2m/p}}.
  \end{equation}
\end{enumerate}
\end{prop}

\begin{proof}
Let $u_m$ and $E(t)$ be given by Lemma~\ref{lemma:ConstructionOfmixingVectorField}, and let $W_{m,J}$ be given by Corollary~\ref{cor:MultipleRotationsInASquare}.
For $x \in \R^2$, set
\begin{equation}\label{eq:CombinedVelocityDefinition}
  b_m(t,x) = 
  \begin{cases}
    u_m(t,x), & 0\leq t \leq1 ;\\
    \displaystyle \frac{\beta 2^{m+J} \omega(2^{-(m+J)})}{\omega(2^{-(m+J)}) + 2^{-2m/p}}  W_{m,J}(x), & 1< t \leq2.
  \end{cases}
\end{equation}
Also set
\begin{equation}\label{eq:CombinedDensityDefinition}
  \rho_m(t,x)=
  \begin{cases}
    \mathbf{1}_{E(t)}(x), & 0\leq t\leq1; \\
    \mathbf{1}_{E_{m,J}}(x), & 1<t\leq2.
  \end{cases}
\end{equation}
Since $u_m$ and $W_{m,J}$ are divergence-free, we deduce that $b_m$ is divergence-free.
Lemma~\ref{lemma:ConstructionOfmixingVectorField}, \eqref{eq:MultipleRotationsNormBounds} and \eqref{eq:MultipleRotationsOsgoodBound} give \eqref{eq:CombinedBlockNormBounds}.
Since $E(1)=E_{m,J}$, the density is continuous in $L^1$ at $t=1$, and it follows from Lemma~\ref{lemma:ConstructionOfmixingVectorField} and Corollary~\ref{cor:MultipleRotationsInASquare} that $\rho_m$ solves the transport equation with $b_m$.
Finally, linearity of the commutator in the vector field and \eqref{eq:MultipleRotationsCommutatorEstimate} give \eqref{eq:OneFullBlockCommutatorEstimate} for $1<t\leq2$.
\end{proof}

\subsection{Proof of Theorem~\ref{thm:CommutatorBlowup}}
\begin{proof}[Proof of Theorem~\ref{thm:CommutatorBlowup}]
Let $c_0>0$ be the minimum of the constant in Lemma~\ref{lemma:ConstructionOfmixingVectorField} and the lower-bound constant in Proposition~\ref{prop:CombinedBuildingBlock}.
Fix a sequence $\beta_n\in(0,1)$ such that
\begin{equation}\label{eq:GluingBetaCondition}
  \sum_{n\geq1}\beta_n\leq 1/4.
\end{equation}
We will choose integers $1\leq m_n<k_n$, with $k_n$ strictly increasing, such that
\begin{equation}\label{eq:GluingParameterConditions}
  2^{-m_n}\leq c_0\beta_n,
  \qquad
  2^{-2m_n/p}G(k_n)\leq c_0\beta_n.
\end{equation}
Set $J_n=k_n-m_n$ and note that, since $G(k_n)-G(m_n)\leq G(k_n)$, Proposition~\ref{prop:CombinedBuildingBlock} with parameters $(\beta_n,m_n,J_n)$ yields a vector field $b_{m_n}$ and a density $\rho_{m_n}$.
By the first inequality of \eqref{eq:GluingParameterConditions} and \eqref{eq:GluingBetaCondition}, there exist points $x_n\in\R^2$ such that the squares $x_n+2Q_{m_n}$ are mutually disjoint and contained in $(-1/4,1/4)^2$.
On $\R^2$, define
\begin{equation}\label{eq:GluedFields}
  b(t,x)=\sum_{n\geq1}b_{m_n}(t,x-x_n),
  \qquad
  \rho(t,x)=\sum_{n\geq1}\rho_{m_n}(t,x-x_n),
\end{equation}
and then periodise these functions to $\T^2=\R^2/\Z^2$.
At any point at most one summand is nonzero, so
\[
  \|b\|_{L^\infty([0,2]\times\T^2)}\leq C\sup_n\beta_n,
  \qquad 0\leq\rho\leq1,
\]
and, using Lemma~\ref{lem:DisjointnessOfSupportsForHajSpaces}, 
\begin{equation}\label{eq:GluedOsgoodNorm}
  \|b\|_{L^\infty([0,2];\dot{M}^{\omega,p}(\T^2))} \leq \sum_{n \geq 1} \left( \|b_{m_n}\|_{L^\infty([0,2];\dot{M}^{\omega,p}(\R^2))} + \frac{\| b_{m_n} \|_{L^\infty([0,2]; L^{p})}}{\omega(2^{-m_n})} \right) \leq C(p) \sum_{n\geq1} \beta_n < \infty,
\end{equation}
that is, $b$ has the desired regularity.
Since each $b_{m_n}$ is divergence-free, $b$ is divergence-free.

The density series converges in $C([0,2];L^1)$, and disjointness gives
\[
  b\rho=\sum_{n\geq1}
  b_{m_n}(t,\cdot-x_n)\rho_{m_n}(t,\cdot-x_n),
\]
and thus $\rho$ solves \eqref{eq:TransportEquation}, with
\[
  \rho(0,\cdot)=\sum_{n\geq1}
  \mathbf{1}_{B_{2^{-(m_n+2)}}(x_n)}.
\]
The initial datum belongs to $\operatorname{BV}(\T^2)$ since the disks have summable perimeters.

Let $\eps_n=\eps_0 2^{-k_n}$, where $\eps_0$ is the constant in Lemma~\ref{lemma:SingleDiskCommutator}.
Since $k_n\geq m_n+1$ and $\eps_0 < 1/5$, we have $\eps_n < 2^{-m_n}/10$.
The translated $n$-th commutator is supported in $x_n + Q_{m_n} + B_{\eps_n}(0)$.
For every $x$ in this neighborhood, the ball $B_{\eps_n}(x)$ lies in $x_n + 2Q_{m_n}$ and meets no other block.
Hence the global commutator agrees there with the translated $n$-th commutator.
By Proposition~\ref{prop:CombinedBuildingBlock}, for every $t\in(1,2]$,
\begin{equation}\label{eq:GluedCommutatorLowerBound}
  \|R_{\eps_n}(b_t,\rho_t)\|_{L^1(\T^2)} \geq \|R_{\eps_n}(b_{m_n}(t,\cdot),\rho_{m_n}(t,\cdot))\|_{L^1(\R^2)}
  \geq c_0\beta_n \frac{2^{k_n-2m_n}\omega(2^{-k_n})}{\omega(2^{-k_n})+2^{-2m_n/p}}.
\end{equation}

It remains to choose $m_n$ and $k_n$.
For each $n$, take $k_n$ sufficiently large and choose $m_n$ so that
\begin{equation}\label{eq:DyadicChoiceOfM}
  \frac12\left(\frac{c_0\beta_n}{G(k_n)}\right)^{p/2}
  <2^{-m_n}\leq
  \left(\frac{c_0\beta_n}{G(k_n)}\right)^{p/2}.
\end{equation}
Since $G(k)\to\infty$ and $G(k)=O(k)$, taking $k_n$ large ensures $1\leq m_n<k_n$ and \eqref{eq:GluingParameterConditions}.
Moreover, $G(k)\omega(2^{-k})\to0$, so we may also require $G(k_n)\omega(2^{-k_n})\leq c_0\beta_n$.
Substituting \eqref{eq:DyadicChoiceOfM} into \eqref{eq:GluedCommutatorLowerBound} then gives, for $t\in(1,2]$,
\[
  \|R_{\eps_n}(b_t,\rho_t)\|_{L^1(\T^2)}
  \geq c\beta_n^p\frac{2^{k_n}\omega(2^{-k_n})}{G(k_n)^{p-1}},
\]
where $c>0$ is independent of $n$.
By \eqref{eq:CommutatorGrowthHypothesis}, we can choose $k_n$ strictly increasing so that the right-hand side is at least $n$.
Then $\eps_n\downarrow0$ and
\[
  \int_0^2\|R_{\eps_n}(b_t,\rho_t)\|_{L^1(\T^2)}\,dt
  \geq\int_1^2\|R_{\eps_n}(b_t,\rho_t)\|_{L^1(\T^2)}\,dt
  \geq n.
\]
\end{proof}

\section{Sharp regularity and mixing}
We start this section by proving the optimal regularity result.
The proof is based on the rescaling and gluing scheme introduced in \cite{AlbertiCrippaMazzucato2019Loss} and refined in \cite[Section~2]{BrueNguyen2021Annalen}. The main difference is that here, the time scales of the individual building blocks are chosen according to the modulus $\omega$, in order to keep the resulting vector field uniformly bounded in $\dot{M}^{\omega,p}$. In addition, we adapt the proof to treat the loss of regularity in the Littlewood-Paley context.

\subsection{Sharpness of $B^{\omega, a}$ regularity}

\begin{proof}[Proof of Theorem~\ref{thm:loss-of-regularity}]
Let $v$ and $\theta$ be given by Theorem~\ref{lipschitz},
with $Q=(-1/2,1/2)^d$, and set $\bar\theta := \theta(0,\cdot)$.
Multiplying $\theta$ by a constant, we may assume that
$\|\bar \theta \|_{L^2(\R^d)}=1$.
For every integer $n \geq 1$, we set
\[
\lambda_n:=2^{-n},
\qquad \gamma_n:=n^{-2},
\qquad \tau_n:=\big( n^2\lambda_n^d \big)^{1/p}.
\]
Choosing $n_0$ sufficiently large, we can find centres $x_n$ such that, for every $n\geq n_0$,
the cubes 
\[
Q_n:=x_n + 4 d \lambda_n Q
\]
are pairwise disjoint and contained in $\T^d$. For every $n\geq n_0$, we also define the function
$F_n:(0,\lambda_n]\to[0,\infty)$ by
\[
F_n(r):=\int_r^{\lambda_n}
\left(1+\frac{\tau_n}{\omega(s)}\right)\,ds.
\]
The function $F_n$ is continuous and strictly decreasing,
with $F_n(\lambda_n)=0$. By the Osgood condition,
$F_n(r)\to+\infty$ as $r\to0^+$, so $F_n$ is a bijection.
For every $t\geq0$, we then define
\begin{equation}\label{eq:loss-time-change}
    \ell_n(t):=F_n^{-1}(t),
\end{equation}
and
\[
s_n(t) := \log\frac{\lambda_n}{\ell_n(t)}.
\]
By construction, we have $\ell_n(0)=\lambda_n$ and $s_n(0)=0$.
Moreover, by differentiating the identity $F_n(\ell_n(t))=t$, we obtain
\begin{equation}\label{eq:ODEForElln}
1=F_n'(\ell_n(t))\ell_n'(t)
=-\left(1+\frac{\tau_n}{\omega(\ell_n(t))}\right)\ell_n'(t) \quad \Rightarrow \quad \ell_n'(t)=-\frac{\omega(\ell_n(t))}{\tau_n+\omega(\ell_n(t))}.
\end{equation}
Hence, and recalling the definition of $s_n$, we have
\[
s_n'(t)=-\frac{\ell_n'(t)}{\ell_n(t)}
=\frac{\omega(\ell_n(t))}
{\ell_n(t)(\tau_n+\omega(\ell_n(t)))}.
\]
In particular, we have
\begin{equation}\label{eq:loss-speed-bounds}
\ell_n(t)s_n'(t) \leq 1,
\qquad s_n'(t) \leq \frac{\omega(\ell_n(t))}{\tau_n\ell_n(t)}.
\end{equation}
For $x \in Q_n$, define
\begin{equation}
b^{(n)}(t,x) := \lambda_n s_n'(t) v\left(s_n(t),\frac{x-x_n}{\lambda_n}\right), \quad
\rho^{(n)}(t,x) := \gamma_n \theta \left(s_n(t),\frac{x-x_n}{\lambda_n}\right),
\end{equation}
and extend both functions by zero outside $Q_n$.
Finally, we define the vector field and the solution
appearing in the statement of the theorem by
\[
b(t,x):=\sum_{n\geq n_0}b^{(n)}(t,x), \qquad \rho(t,x):=\sum_{n\geq n_0}\rho^{(n)}(t,x).
\]
By \eqref{eq:DecayOfSmoothBuildingBlock} and
\eqref{eq:loss-speed-bounds}, we have
\[
\|b^{(n)}(t)\|_{L^\infty} \lesssim  \lambda_n s_n'(t)e^{-s_n(t)} = \ell_n(t)s_n'(t)\leq 1.
\]
Since the supports are pairwise disjoint, it follows that $\|b\|_{L^\infty([0,\infty)\times\T^d)} \lesssim 1$. 
By using \eqref{eq:InterpolationEstimateAtScale}, for $r=\ell_n(t)$, we have
\[
\| b^{(n)} \|_{ \dot{M}^{\omega,p} } \lesssim_{d,p} \frac{\max\{\|b^{(n)}\|_{L^p}, \ell_n(t)\| \nabla b^{(n)} \|_{L^p}\} }{\omega(\ell_n(t))}.
\]
By definition of $b_n$ and \eqref{eq:LipRegOfSmoothBuildingBlock}
and \eqref{eq:DecayOfSmoothBuildingBlock}, we also have
\[
\begin{aligned}
\|b^{(n)}(t)\|_{L^p} &\lesssim s_n'(t)\lambda_n^{1+d/p}e^{-s_n(t)} = s_n'(t)\ell_n(t)\lambda_n^{d/p},\\
\|\nabla b^{(n)}(t)\|_{L^p}
&\lesssim s_n'(t)\lambda_n^{d/p}.
\end{aligned}
\]
By combining the previous inequalities and using \eqref{eq:loss-speed-bounds}, we get
\[
\|b^{(n)}\|_{\dot{M}^{\omega,p}} \lesssim_{d,p} s_n'(t)\frac{\ell_n(t)\lambda_n^{d/p}}{\omega(\ell_n(t))}
\leq \frac{\lambda_n^{d/p}}{\tau_n} = n^{-2/p}.
\]
Then, since the supports of $b_n$ are well separated, using Lemma~\ref{lem:DisjointnessOfSupportsForHajSpaces}, we deduce that \begin{equation}
    \| b \|_{ \dot{M}^{\omega,p} }^p \leq \sum_{n \geq n_0} \left( \| b^{(n)} \|_{ \dot{M}^{\omega,p} } + \frac{\|b^{(n)}\|_p}{\omega(d \lambda_n)} \right)^p \lesssim_{d,p} 1.
\end{equation}
Thus $b\in L^\infty([0,\infty);\dot{M}^{\omega,p}(\T^d;\R^d))$. For $\rho$, since the supports of $\rho^{(n)}$ are disjoint, we have
\begin{equation}
\| \rho \|_{L^\infty([0,\infty)\times\T^d)} \leq \sup_{n\geq n_0} \gamma_n  \|\bar \theta \|_{L^\infty(\R^d)}.
\end{equation}
Each $b_n$ is divergence-free and each pair $(b^{(n)},\rho^{(n)})$ solves
the transport equation. Since the supports of the pairs $(b^{(n)},\rho^{(n)})$ are disjoint, we have that
$b\rho=\sum_{n\geq n_0} b^{(n)} \rho^{(n)}$, $b$ is divergence-free and $\rho$
is a bounded weak solution of \eqref{eq:TransportEquation}, with
initial datum
\begin{equation}
\rho^{(0)} (x) = \sum_{n\geq n_0} \gamma_n \bar \theta \left( \frac{x-x_n}{\lambda_n} \right).
\end{equation} 
The same argument also gives $\rho\in C([0,\infty);L^2(\T^d))$.
Mutual disjointness of supports gives
\[
\|\nabla \rho^{(0)} \|_{L^d}^d
=
\|\nabla \bar \theta \|_{L^d(\R^d)}^d
\sum_{n\geq n_0}\gamma_n^d
<\infty.
\]
Hence $\rho_0 \in L^\infty(\T^d)\cap W^{1,d}(\T^d)$.
Finally, we fix $t>0$ and $a>p/2$ and show $\rho(t,\cdot)\notin B^{\omega,a}(\mathbb{T}^d)$.
By definition of $\rho_n$ and \eqref{eq:DecayOfSmoothBuildingBlock}, we get
\begin{equation}\label{eq:loss-block-norms}
\|\rho^{(n)}(t)\|_{L^2} = \gamma_n\lambda_n^{d/2}, \quad
\|\rho^{(n)}(t)\|_{\dot H^{-1}} \lesssim \gamma_n \lambda_n^{1 + d/2} \e^{- s_n(t)}=\gamma_n \lambda_n^{d/2} \ell_n(t).
\end{equation}
Using the disjointness of the supports,
\begin{equation}
    \int_{\T^d} \rho(t) \rho^{(n)}(t) \, dx = \|\rho^{(n)}(t) \|_{L^2}^2 = \gamma_n^2\lambda_n^d.
\end{equation}
Let $k_n \geq 1$ be an integer depending on $n$ to be selected and use the decomposition $\rho = \rho_{k_{n}-1}^{\leq} + \rho_{k_{n}}^{\geq}$ to deduce
\begin{equation}
\begin{aligned}
\gamma_n^2\lambda_n^d &\leq \|\rho^{(n)}(t)\|_{L^2} \|\rho_{k_n}^{\geq}(t)\|_{L^2} + \|\nabla \rho_{k_n-1}^{\leq}(t)\|_{L^2} \|\rho^{(n)}(t)\|_{\dot H^{-1}}\\
&\leq \gamma_n\lambda_n^{d/2}\|\rho_{k_n}^{\geq}(t)\|_{L^2} + C \|\rho(t)\|_{L^2}2^{k_n} \gamma_n \lambda_n^{d/2}\ell_n(t).
\end{aligned}
\end{equation}
Taking $k_n := \lfloor \log_2(1/\ell_n(t))\rfloor - d n$, the last term is bounded by $\frac{1}{2} \gamma_n^2 \lambda_n^d$ and hence
\begin{equation}
    \gamma_n^2 \lambda_n^d \leq 2 \| \rho_{k_n}^{\geq} \|_{L^2} \gamma_n \lambda_n^{d / 2} \quad \Rightarrow \quad \frac{1}{2} \gamma_n \lambda_n^{d / 2} \leq \| \rho_{k_n}^{\geq} \|_{L^2}.
\end{equation}
Then
\begin{equation}
    \| \rho(t) \|^2_{B^{\omega, a}(\T^d)} \gtrsim G(k_n)^{2 a} \| \rho_{k_n}^{\geq}(t) \|^2_{L^2(\T^d)} \gtrsim G(k_n)^{2 a} \gamma_n^2 \lambda_n^d.
\end{equation}
To estimate $G(k_n)^{2 a}$ from below, we first observe that
\begin{equation}
    \dfrac{d}{dt}\left[ G\Big(\log_2(1/t)\Big) \right] = - \frac{1}{\omega(t)} \quad \Rightarrow \quad \dfrac{d}{dt}\left[ G\Big(\log_2(1/\ell_n(t))\Big) \right] = - \frac{\ell_n^\prime(t)}{\omega(\ell_n(t))} \stackrel{\eqref{eq:ODEForElln}}{=} \frac{1 + \ell_n^\prime(t)}{\tau_n}
\end{equation}
so that by integrating from $0$ to $t$ and using $\ell_n(0) = \lambda_n = 2^{-n}$, we get
\begin{equation}
    G\Big(\log_2(1/\ell_n(t))\Big) - G(n) = \frac{t + \ell_n(t) - \lambda_n}{\tau_n} \geq \frac{t - \lambda_n}{\tau_n}.
\end{equation}
Then, by Lemma~\ref{lem}, 
\begin{equation}
    G(k_n) \geq G\Big(\log_2(1/\ell_n(t))\Big) - (dn+1) G(1) \geq \frac{t}{2 \tau_n} \quad \forall t \geq C( \lambda_n + n\tau_n)
\end{equation}
so that 
\begin{equation}
\begin{aligned}
\|\rho(t)\|_{B^{\omega,a}}^2 &\gtrsim t^{2a}\gamma_n^2\lambda_n^d\tau_n^{-2a} \gtrsim t^{2a}n^{-4-\frac{4a}{p}}
2^{d n \left(\frac{2a}{p}-1\right)} \quad \forall t \geq C( \lambda_n + n\tau_n).
\end{aligned}
\end{equation}
Since $a > p / 2$, letting $n \to \infty$ gives $\rho(t, \cdot) \not \in B^{\omega, a}(\T^d)$.
\end{proof}

\subsection{Lower bounds on mixing efficiency}

\begin{lem}\label{lem:interpolation}
Let $d \geq 2$, $a \in (0,\infty)$ and let $\omega$ be an Osgood modulus of continuity.
Let $f\in L^2(\mathbb T^d)\cap B^{\omega,a}(\mathbb T^d)$ be nonzero and mean zero.
Then
\begin{equation}\label{eq:interpolation}
  \|f\|_{L^2}^2 \leq C(d,a,\omega) \left(  \|f\|_{\dot H^{-1}}^2 + \|f\|_{B^{\omega,a}}^2 G \left( \log_2 \left( 2+ \frac{\|f\|_{B^{\omega,a}}}{\|f\|_{\dot H^{-1}}} \right) \right)^{-2a} \right).
\end{equation}
\end{lem}

\begin{proof}
For every real $\lambda>0$, split the dyadic sum of \eqref{eq:PlancheralEquivalence} into $k\leq\lambda$ and $k>\lambda$.
The low-frequency part is controlled by the $\dot H^{-1}$ norm, while $G(k)\geq G(\lambda)$ on the high-frequency part.
The Littlewood--Paley estimates give, uniformly in $\lambda$,
\begin{equation}\label{eq:split}
\begin{split}
  \|f\|_{L^2}^2 &\lesssim \sum_{k \leq \lambda}2^{2(\lambda - k)} \| f_k \|_{L^2}^2 + G(\lambda)^{-2a} \sum_{k > \lambda} G(k)^{2a} \| f_k \|_{L^2}^2 \\
  &\lesssim 2^{2\lambda}\|f\|_{\dot H^{-1}}^2 + G(\lambda)^{-2a}\|f\|_{B^{\omega,a}}^2.
\end{split}
\end{equation}
Set
\[
  \lambda:=\frac{1}{2} \log_2 \left(2+ \frac{\|f\|_{B^{\omega,a}}}{\|f\|_{\dot H^{-1}}}\right)\geq\frac12.
\]
Since
\[
  2^{2\lambda}\|f\|_{\dot H^{-1}}^2 = 2\|f\|_{\dot H^{-1}}^2 + \|f\|_{B^{\omega,a}}\|f\|_{\dot H^{-1}}
\]
and $G(2\lambda) \leq 2 G(\lambda)$, \eqref{eq:split} gives
\[
  \|f\|_{L^2}^2\lesssim\|f\|_{\dot H^{-1}}^2+\|f\|_{B^{\omega,a}}\|f\|_{\dot H^{-1}}+G(2\lambda)^{-2a}\|f\|_{B^{\omega,a}}^2.
\]
To obtain \eqref{eq:interpolation}, we show that the second term can be bounded by the sum of the first and the third.
First, applying the elementary inequality $r \leq 1 + \frac{r^2}{1+r}$ for $r \geq 0$ with $r = \|f\|_{B^{\omega,a}} / \|f\|_{\dot H^{-1}}$, and multiplying by $\|f\|_{\dot H^{-1}}^2$, gives
\begin{equation}\label{eq:EstimatingMiddleQuantityIntermediateStep}
  \| f \|_{B^{\omega,a}} \| f \|_{\dot H^{-1}} \leq \| f \|_{\dot H^{-1}}^2 + \frac{\|f\|_{B^{\omega,a}}^2}{1 + \| f \|_{B^{\omega,a}}/\| f \|_{\dot H^{-1}}}.
\end{equation}
Thus, we only have to treat the last term.
Using that $\omega$ is concave and $\omega(0) = 0$, it follows that $( G(\log_2(2+r)) )^{2 a} \lesssim_{a,\omega} 1 + r$.
Applying this with $r = \| f \|_{B^{\omega,a}}/\| f \|_{\dot H^{-1}}$ and recalling the choice of $\lambda$ implies that the last term of \eqref{eq:EstimatingMiddleQuantityIntermediateStep} is estimated by $G(2\lambda)^{-2a}\| f \|_{B^{\omega,a}}^2$.
\end{proof}

\begin{cor}\label{cor:lower}
Let $d \geq 2$, $1 < p < \infty$, $a \in (0,\infty)$ and let $\omega$ be an Osgood modulus of continuity.
For every nonzero, mean-zero $f\in L^2(\mathbb T^d)\cap B^{\omega,a}(\mathbb T^d)$,
\begin{equation}\label{eq:direct-lower}
  \|f\|_{\dot H^{-1}} \geq \|f\|_{L^2} \exp \left[ - C  G^{-1} \left( C(a) \left( \frac{\|f\|_{B^{\omega,a}}}{\|f\|_{L^2}}\right)^{1/a} \right) \right].
\end{equation}
\end{cor}

\begin{proof}
If $\|f\|_{\dot H^{-1}}$ is larger than a small fixed multiple of $\|f\|_{L^2}$, the conclusion is immediate, since the exponential factor in \eqref{eq:direct-lower} is at most $1$.
Otherwise, the first term on the right-hand side of \eqref{eq:interpolation} can be absorbed into the left-hand side, giving
\[
  G \left( \log_2  \left( 2+
    \frac{\|f\|_{B^{\omega,a}}}{\|f\|_{\dot H^{-1}}} \right) \right) \lesssim_a \left(\frac{\|f\|_{B^{\omega,a}}}{\|f\|_{L^2}}\right)^{1/a}.
\]
Since the logarithm is at least $1$ and $G(1)>0$, this also implies $\|f\|_{B^{\omega,a}} \gtrsim_a \|f\|_{L^2}$.
Applying $G^{-1}$ and taking the exponential yields
\[
  \|f\|_{\dot H^{-1}} \geq \|f\|_{B^{\omega,a}} \exp \left[-C G^{-1} \left( C(a)\left(\frac{\|f\|_{B^{\omega,a}}}{\|f\|_{L^2}}\right)^{1/a} \right)\right],
\]
which proves \eqref{eq:direct-lower}.
\end{proof}

\begin{proof}[Proof of Theorem~\ref{thm:Mixing}]
    By \cite[Theorem 1.2]{dephilippis_franchi_2026_osgood}, for any $a \in [\frac{1}{2}, \frac{p}{2}]$,
    \begin{equation}
        \| \rho(t, \cdot) \|_{B^{\omega, a}} \lesssim_{d,\omega,p,a} \left( \int_0^t \| b(s, \cdot) \|_{\dot{M}^{\omega, p}} \, ds \right)^{a} \| \rho_0 \|_{L^\infty} + \| \rho_0 \|_{B^{\omega, a}}
    \end{equation}
    Dividing by $\| \rho(t, \cdot) \|_{L^2}$ and using the elementary inequality $x^{a} + y^{a} \lesssim_a (x+y)^{a}$, we get
\begin{equation}
    \dfrac{\| \rho(t, \cdot) \|_{B^{\omega, a}}}{\| \rho(t, \cdot) \|_{L^2}} \lesssim_{d,\omega,p,a}
    \left( \left( \int_0^t \| b(s, \cdot) \|_{\dot{M}^{\omega, p}} \, ds \right)
    \left( \frac{\| \rho_0 \|_{L^\infty}}{\| \rho_0 \|_{L^2}} \right)^{1/a}
    + \left( \frac{\| \rho_0 \|_{B^{\omega,a}}}{\| \rho_0 \|_{L^2}} \right)^{1/a} \right)^a.
\end{equation}
Taking $a=\frac{1+p}{4}$, and setting
\[
A_\infty := \left(\frac{\|\rho_0\|_{L^\infty}}{\|\rho_0\|_{L^2}}\right)^{\frac{4}{1+p}}, \qquad
A_{\omega,p} := \left(\frac{\|\rho_0\|_{B^{\omega,(1+p)/4}}}{\|\rho_0\|_{L^2}}\right)^{\frac{4}{1+p}},
\]
Corollary~\ref{cor:lower} gives
\begin{equation}
    \frac{\|\rho(t,\cdot)\|_{\dot H^{-1}}}
    {\|\rho(t,\cdot)\|_{L^2}}
    \geq
    \exp\left[
    -C G^{-1}\left(
    C\left(
    A_\infty \int_0^t
    \|b(s,\cdot)\|_{\dot M^{\omega,p}}\,ds
    + A_{\omega,p}
    \right)\right)\right].
\end{equation}
Since
\[
G'(r)=(\log 2)\frac{2^{-r}}{\omega(2^{-r})}
\]
is positive and nonincreasing by concavity of $\omega$,
$G^{-1}$ is increasing and convex. In particular,
\[
G^{-1}(x+y)\leq G^{-1}(2x)+G^{-1}(2y),
\qquad x,y\geq 0.
\]
Moreover, convexity and $G^{-1}(0)=0$ imply
$C G^{-1}(z)\leq G^{-1}(Cz)$ for $C\geq 1$.
Thus, up to enlarging $C$,
\[
\frac{\|\rho(t,\cdot)\|_{\dot H^{-1}}}
{\|\rho_0\|_{L^2}}
\geq
\Phi_\omega(C A_{\omega,p})
\Phi_\omega\left(
C A_\infty \int_0^t
\|b(s,\cdot)\|_{\dot M^{\omega,p}} \, ds
\right),
\]
which concludes the proof.
\end{proof}

\begin{rem}[Sharpness of the mixing rate]\label{rem:SharpMixing}

The mixing rate in Theorem~\ref{thm:Mixing} is essentially sharp, as we show by adapting the construction in the proof of Theorem~\ref{thm:loss-of-regularity} to a single building block. Let $v$ and $\theta$ be given by Theorem~\ref{lipschitz}, and set
\[
    F(r)=\int_r^1\left(1+\frac{1}{\omega(s)}\right)\,ds, \quad\ell(t)=F^{-1}(t),\quad s(t)=\log(1/\ell(t)).
\]
Then, $\ell(t) \leq \Phi_\omega(ct)$. We now define
\[
    b(t,x):=s^\prime(t)v(s(t),x), \quad \rho(t,x):=\theta(s(t),x)
\]
By \eqref{eq:DecayOfSmoothBuildingBlock}, we obtain $\|\rho(t)\|_{\dot H^{-1}}\lesssim\ell(t)\leq\Phi_\omega(ct)$.
\end{rem}

\appendix

\section{Regularity of the scalar}\label{sec:RegOfScalar}
In this appendix, we establish Haj\l{}asz-type estimates for solutions to the transport equation with Osgood-Sobolev vector fields. 
We give two proofs. 
The Lagrangian proof uses the flow representation of solutions, justified by the well-posedness result in \cite[Theorem 1.1]{dephilippis_franchi_2026_osgood}. 
The Eulerian proof uses Littlewood-Paley theory and the propagation of regularity established in \cite[Theorem 1.2]{dephilippis_franchi_2026_osgood}. 
Although it does not give the sharp estimate for $p>2$, we include it as a complementary approach.

\subsection{Lagrangian approach}
Here we provide a Lagrangian argument to show that solutions $\rho$ to \eqref{eq:TransportEquation} with $b \in L^\infty([0,2]; \dot{M}^{\omega,p}(\T^d;\R^d)\cap L^\infty(\T^d;\R^d) )$ divergence-free and with sufficiently regular initial datum (including $\operatorname{BV}$) satisfy
\begin{equation}\label{eq:LagrangianWhatRegulartityToProve}
 |\rho(t,x) - \rho(t,y)| \leq \frac{q(t,x) + q(t,y)}{G(- \log_2 (\mathsf{d}(x,y)\wedge\tfrac12))^{p-1}} \text{ with $q(t,\cdot) \in L^{p^\prime}(\T^d)$}.
\end{equation}
Solutions $\rho$ are given by $\rho(t,x) = \rho_0(X_t^{-1}(x))$, where $X_t$ is the regular Lagrangian flow of $b$. It satisfies
\begin{equation}
    \left\{
    \begin{array}{ll}
        \dfrac{d}{dt} X_t(x) = b(t, X_t(x)); \smallskip \\
        X_0(x) = x; \smallskip \\
        (X_t)_{\#} \mathcal{L}^d = \mathcal{L}^d
    \end{array}
    \right.
\end{equation}
If $\rho_0$ is sufficiently regular, then it suffices to prove \eqref{eq:LagrangianWhatRegulartityToProve} for $X_t^{-1}$. Up to reversing the vector field in time, this is equivalent to proving \eqref{eq:LagrangianWhatRegulartityToProve} for $X_t$.
With $b \in L^1_t \dot{M}_x^{\omega,p}$, we have the differential inequality
\begin{equation}
    \dfrac{d}{dt} \mathsf{d}(X_t(x), X_t(y)) \leq \big( g_t(X_t(x)) + g_t(X_t(y)) \big) \omega \big( \mathsf{d}(X_t(x), X_t(y)) \big).
\end{equation}
Then, the Bihari-LaSalle inequality gives, whenever $X_t(x) \neq X_t(y)$,
\begin{equation}\label{eq:ConsequenceOfApplyingBihariLaSalle}
    G(-\log_2 \mathsf{d}(X_t(x), X_t(y)) ) \geq G(- \log_2 \mathsf{d}(x,y) ) - h_t(x) - h_t(y)
\end{equation}
with
\begin{equation}
    h_t(x) = \int_0^t g_s (X_s(x)) \, ds.
\end{equation}
Multiplying \eqref{eq:ConsequenceOfApplyingBihariLaSalle} by $\mathsf{d}(X_t(x), X_t(y))^\beta$ for some $\beta>0$ to be chosen and using that \\
$\sup_{0<s<1} s^\beta G(- \log_2 (s) ) \lesssim_{\beta,\omega} 1$, we get
\begin{equation}
    C(\beta) \geq \Big( G(- \log_2 \mathsf{d}(x,y) ) - h_t(x) - h_t(y) \Big) \mathsf{d}(X_t(x) , X_t(y))^\beta.
\end{equation}
Using that $\mathsf{d}(X_t(x), X_t(y)) \leq C(d)$, we get
\begin{equation}
    \mathsf{d}(X_t(x), X_t(y))^\beta G(- \log_2 \mathsf{d}(x,y) ) \leq C(d,\beta) ( 1 + h_t(x) + h_t(y) ).
\end{equation}
Thus, dividing by $G(-\log_2 \mathsf{d}(x,y))$ and raising both sides to the power $1/\beta$, we obtain
\begin{equation}
    \mathsf{d}(X_t(x), X_t(y)) \leq \dfrac{C(d,\beta) ( 1 + h_t(x) + h_t(y) )^{1/\beta}}{G(- \log_2 \mathsf{d}(x,y) )^{1/\beta}}.
\end{equation}
Choosing $\beta = \frac{1}{p-1}$ and $q(t,x) = C(d, \beta) (1 + h_t(x)^{p-1})$, we get
\begin{equation}
    \mathsf{d}(X_t(x), X_t(y)) \leq \dfrac{q(t,x) + q(t,y)}{G(- \log_2 \mathsf{d}(x,y) )^{p-1}}
\end{equation}
as desired.

\subsection{Littlewood-Paley approach}

We begin by recalling some properties of the function $G$,
which are essentially contained in
\cite{dephilippis_franchi_2026_osgood}.
\begin{lem}\label{lem}
Let $\omega$ be an Osgood modulus of continuity and let
\[
 G(x) := \int_{2^{-x}}^1\frac{ds}{\omega(s)}, \qquad x \geq 0.
\]
Then $G$ is increasing, $G(k) \to +\infty$ as $k \to \infty$, and, for every $k\geq1$,
\[
 G(k)\leq G(k+1)\leq G(k)+G(1),
\]
and
\[
 G(k+1)-G(k)\sim \frac{2^{-k}}{\omega(2^{-k})}.
\]
In particular, for every fixed $L \in \mathbb N$,
\[
 G(k+L)\sim_L G(k).
\]
Moreover, for every $\alpha > 0$ and every $1 \leq m \leq N$,
\[
 \sum_{k=m}^{N} \frac{2^{-k}}{\omega(2^{-k})}G(k)^{\alpha-1} \sim_\alpha G(N+1)^\alpha-G(m)^\alpha.
\]
\end{lem}

\begin{proof}
The first two estimates follow from \cite[Lemmas 2.3 and 2.5]{dephilippis_franchi_2026_osgood}. The estimate for $G(k+L)$ follows by iteration. Finally, since $G(k+1) \sim G(k)$, the mean value theorem gives
\[
G(k+1)^\alpha - G(k)^\alpha \sim_\alpha G(k)^{\alpha-1}\big(G(k+1)-G(k)\big) \sim_\alpha G(k)^{\alpha-1}\frac{2^{-k}}{\omega(2^{-k})},
\]
and summing over $k=m , \ldots , N$ yields the last claim.
\end{proof}
We now prove the pointwise consequences of the Littlewood-Paley
regularity. We recall that the corresponding norm is defined by
\[
\| f \|_{B^{\omega, a}}^2 := \sum_{k = 0}^\infty G(k)^{2 a} \| f_k \|_{L^2}^2.
\]
Before proving the inequality, we need an interpolation lemma.
\begin{lem}\label{interpolationG}
      Let $a>0$, $b\geq 0$ and $r\geq 2$ be given, and let $\rho \in L^\infty(\mathbb{T}^d) \cap B^{\omega,a}$.

\begin{enumerate}
    \item[a)] If $2a = br$, then it holds
    \begin{equation}\label{eq:interpolation_a}
     \|\sup_{k \geq 1} G(k)^{b}  |\rho_k^{\geq}| \|_{L^r} 
    \lesssim \| \rho \|_{L^\infty}^{1 - \frac{b}{a}} \bw^{\frac{b}{a}}.
    \end{equation}

    \item[b)]If $2a = br$, then it holds
    \begin{equation}\label{eq:interpolation_b}
    \|\sup_{k \geq 0} G(k)^{b} |\rho_k| \|_{L^r} 
    \lesssim \| \rho \|_{L^\infty}^{1 - \frac{b}{a}} \bw^{\frac{b}{a}}.
    \end{equation}
\end{enumerate}

\end{lem}
\begin{proof}
The proof follows the same argument as \cite[Lemma~2]{meyer2024propagation} for the case $\omega(r)=r$, combined with Lemma~\ref{lem}.
\end{proof}

We now prove the pointwise inequality.
In the range $1<p<2$, we will also need to interpolate with the $L^\infty$ bound to obtain a pointwise bound involving a function in $L^{p'}$. 

\begin{lem}\label{lem:upointwise}
Let $a>0$, $1<p<\infty$, and let $\rho\in L^\infty(\T^d)\cap B^{\omega,a}(\T^d)$.

\begin{enumerate}
    \item[a)] If $p\geq2$, there exists a nonnegative function
    $q\in L^2(\T^d)$ such that
    \begin{equation}\label{eq:upointwise-L2}
    |\rho(x)-\rho(y)| \lesssim \frac{q(x)+q(y)} {G \big(-\log_2( \mathsf{d}(x,y) \wedge \frac{1}{2})\big)^a}
    \end{equation}
    for almost every $x, y \in \T^d$, with $x \neq y$, and
    \[
    \|q\|_{L^2} \lesssim \|\rho\|_{B^{\omega,a}}.
    \]
    \item[b)] If $1<p\leq2$, there exists a nonnegative function
    $q\in L^{p'}(\T^d)$ such that
    \begin{equation}\label{eq:upointwise-Lpprime}
    |\rho(x)-\rho(y)| \lesssim \frac{q(x)+q(y)}{G\big(-\log_2( \mathsf{d}(x,y) \wedge\tfrac12)\big)^{2a/p'}}
    \end{equation}
    for almost every $x,y\in\T^d$, with $x \neq y$, and
    $\|q\|_{L^{p'}} \lesssim \|\rho\|_{L^\infty}^{1-\frac{2}{p'}} \|\rho\|_{B^{\omega,a}}^{\frac{2}{p'}}$.
\end{enumerate}
\end{lem}
\begin{proof}
We prove both statements simultaneously. Set
\[
(r,\beta) =
\begin{cases}
(2,a), & p \geq 2,\\
(p',2a/p'), & 1 < p \leq 2.
\end{cases}
\]
Note that in both cases we have $r \geq 2$ and $\beta r = 2 a$.
We define
\[
q(x) := \sup_{k\geq1}G(k)^\beta |\rho_k^{\geq}(x)| + M \left( \sup_{k\geq1}G(k)^\beta |\rho_k| \right)(x).
\]
By Lemma~\ref{interpolationG}, together with the $L^r$ boundedness of
the Hardy--Littlewood maximal operator, we obtain
\begin{align*}
\|q\|_{L^r} &\lesssim \left\| \sup_{k\geq1}G(k)^\beta |\rho_k^{\geq}| \right\|_{L^r} + \left\| \sup_{k\geq1}G(k)^\beta |\rho_k| \right\|_{L^r}\\
&\lesssim \|\rho\|_{L^\infty}^{1-\frac{\beta}{a}} \|\rho\|_{B^{\omega,a}}^{\frac{\beta}{a}}\\
&= \|\rho\|_{L^\infty}^{1-\frac{2}{r}} \|\rho\|_{B^{\omega,a}}^{\frac{2}{r}}.
\end{align*}
We can assume that $0 < |x-y| \leq \frac{1}{2}$, and choose $J\geq1$ such that
\[
2^{-J-1}<|x-y|\leq2^{-J}.
\]
We decompose
\[
\rho=\rho_{J+1}^{\leq}+\rho_{J+2}^{\geq}.
\]
For the high-frequency part, by the definition of $q$ and the
monotonicity of $G$, we have
\begin{align*}
|\rho_{J+2}^{\geq}(x)-\rho_{J+2}^{\geq}(y)| &\leq |\rho_{J+2}^{\geq}(x)|+|\rho_{J+2}^{\geq}(y)|\\
&\leq \frac{q(x)+q(y)}{G(J+2)^\beta}\\
&\leq \frac{q(x)+q(y)}{G(J)^\beta}.
\end{align*}
It remains to estimate the low-frequency part. By Peetre's maximal inequality \cite[Lemma~14.6.1]{hytonenvanneervenveraarweis2023}, we have, for every $j\geq1$,
\begin{equation}\label{eq:peetre-gradient}
|\nabla \rho_j(x-z)| \lesssim 2^j(1+2^j|z|)^N M(|\rho_j|)(x),
\end{equation}
where $N>d$ is fixed.
For $j\leq J+1$ and $\theta\in[0,1]$, $2^j\theta|x-y|\leq2$, and therefore
\[
|\nabla \rho_j(x+\theta(y-x))| \lesssim 2^jM(|\rho_j|)(x).
\]
Finally, by using the fundamental theorem of calculus, we have
\begin{align*}
|\rho_{J+1}^{\leq}(x)-\rho_{J+1}^{\leq}(y)|
&\lesssim \mathsf{d}(x,y) \sum_{j=1}^{J+1}2^jM(|\rho_j|)(x)\\
&\leq \mathsf{d}(x,y) \,q(x) \sum_{j=1}^{J+1}\frac{2^j}{G(j)^\beta}\\
&\lesssim \frac{q(x)}{G(J)^\beta},
\end{align*}
where in the last inequality we used Lemma~\ref{lem} and $|x-y|\sim 2^{-J}$.
Since $2^{-J-1}<|x-y|\leq 2^{-J}$, Lemma~\ref{lem} gives
\[
G(J) \sim G\big(-\log_2 \mathsf{d}(x,y) \big),
\]
and therefore
\[
|\rho(x)-\rho(y)| \lesssim \frac{q(x)+q(y)}{G\bigl(-\log_2\mathsf{d}(x,y) \bigr)^\beta}.
\]
\end{proof}

\printbibliography
 
\end{document}